\input amstex

\documentstyle{amsppt}
  \magnification=1100
  \hsize=6.2truein
  \vsize=9.0truein
  \hoffset 0.1truein
  \parindent=2em

\NoBlackBoxes


\font\eusm=eusm10                   


\font\eusms=eusm7                       

\font\eusmss=eusm5                      


\newcount\theTime
\newcount\theHour
\newcount\theMinute
\newcount\theMinuteTens
\newcount\theScratch
\theTime=\number\time
\theHour=\theTime
\divide\theHour by 60
\theScratch=\theHour
\multiply\theScratch by 60
\theMinute=\theTime
\advance\theMinute by -\theScratch
\theMinuteTens=\theMinute
\divide\theMinuteTens by 10
\theScratch=\theMinuteTens
\multiply\theScratch by 10
\advance\theMinute by -\theScratch

\def\today{{\number\day\space
 \ifcase\month\or
  January\or February\or March\or April\or May\or June\or
  July\or August\or September\or October\or November\or December\fi
 \space\number\year}}

\define\Ah{{\widehat A}}

\define\Ao{{A\oup}}

\define\biggnm#1{
  \bigg|\bigg|#1\bigg|\bigg|}

\define\bignm#1{
  \big|\big|#1\big|\big|}

\define\clspan{\overline\lspan}

\define\Cpx{\bold C}

\define\eqdef{{\;\overset\text{def}\to=\;}}

\define\Eto#1{E_{(\to{#1})}}

\define\fpamalg#1{{\dsize\;\operatornamewithlimits*_{#1}\;}}

\define\fpiamalg#1{{\tsize\;({*_{#1}})_{\raise-.5ex\hbox{$\ssize\iota\in I$}}}}

\define\freeprod#1#2{\mathchoice
     {\operatornamewithlimits{\ast}_{#1}^{#2}}
     {\raise.5ex\hbox{$\dsize\operatornamewithlimits{\ast}
      _{#1}^{#2}$}\,}
     {\text{oops!}}{\text{oops!}}}

\define\freeprodi{\mathchoice
     {\operatornamewithlimits{\ast}
      _{\iota\in I}}
     {\raise.5ex\hbox{$\dsize\operatornamewithlimits{\ast}
      _{\sssize\iota\in I}$}\,}
     {\text{oops!}}{\text{oops!}}}

\define\freeprodvni{\mathchoice
      {\operatornamewithlimits{\overline{\ast}}
       _{\iota\in I}}
      {\raise.5ex\hbox{$\dsize\operatornamewithlimits{\overline{\ast}}
       _{\sssize\iota\in I}$}\,}
      {\text{oops!}}{\text{oops!}}}

\define\GNS{{\text{\rm GNS}}}

\define\Hil{{\mathchoice
     {\text{\eusm H}}
     {\text{\eusm H}}
     {\text{\eusms H}}
     {\text{\eusmss H}}}}


\define\Hilto#1{\Hil_{(\to{#1})}}

\define\Ic{{\Cal I}}

\define\id{\text{\rm id}}

\define\Keu{{\KHil}}

\define\KHil{{\mathchoice
     {\text{\eusm K}}
     {\text{\eusm K}}
     {\text{\eusms K}}
     {\text{\eusmss K}}}}

\define\ld#1{{\hbox{..}(#1)\hbox{..}}}

\define\Leu{{\mathchoice
     {\text{\eusm L}}
     {\text{\eusm L}}
     {\text{\eusms L}}
     {\text{\eusmss L}}}}

\define\lrnm#1{\left\|#1\right\|}

\define\lspan{\text{\rm span}@,@,@,}

\define\nm#1{\|#1\|}

\define\Nats{\Naturals}

\define\Naturals{{\bold N}}

\define\otdt{\otimes\cdots\otimes}

\define\otdts#1{\otimes_{#1}\cdots\otimes_{#1}}

\define\otimesdd{{\,{\overset{..}\to\otimes}}\,}

\define\oup{^{\text{\rm o}}}

\define\owedge{{
     \operatorname{\raise.5ex\hbox{\text{$
     \ssize{\,\bigcirc\llap{$\ssize\wedge\,$}\,}$}}}}}

\define\owedgeo#1{{
     \underset{\raise.5ex\hbox
     {\text{$\ssize#1$}}}\to\owedge}}

\define\Po{P\oup}

\define\Pto#1{{P_{(\to{#1})}}}


\define\pup#1#2{{{\vphantom{#2}}^{#1}\!{#2}}\vphantom{#2}}

\define\QED{\newline
            \line{$\hfill$\qed}\enddemo}

\define\red{{\text{\rm red}}}

\define\restrict{\lower .3ex
     \hbox{\text{$|$}}}

\define\smd#1#2{\underset{#2}\to{#1}}

\define\smdb#1#2{\undersetbrace{#2}\to{#1}}

\define\smdbp#1#2#3{\overset{#3}\to
     {\smd{#1}{#2}}}

\define\smdbpb#1#2#3{\oversetbrace{#3}\to
     {\smdb{#1}{#2}}}

\define\smdp#1#2#3{\overset{#3}\to
     {\smd{#1}{#2}}}

\define\smdpb#1#2#3{\oversetbrace{#3}\to
     {\smd{#1}{#2}}}

\define\smp#1#2{\overset{#2}\to
     {#1}}

\define\st{{\tilde s}}

\define\Tcirc{\bold T}

\define\tocdots
  {\leaders\hbox to 1em{\hss.\hss}\hfill}


  \newcount\mycitestyle \mycitestyle=1 

  \newcount\bibno \bibno=0
  \def\newbib#1{\advance\bibno by 1 \edef#1{\number\bibno}}
  \ifnum\mycitestyle=1 \def\cite#1{{\rm[\bf #1\rm]}} \fi
  \def\scite#1#2{{\rm[\bf #1\rm, #2]}}


  \newcount\ignorsec \ignorsec=0
  \def\notasec{\ignorsec=1}

  \newcount\secno \secno=0
  \def\newsec#1{\procno=0 \subsecno=0 \ignorsec=0
    \advance\secno by 1 \edef#1{\number\secno}
    \edef\currentsec{\number\secno}}

  \newcount\subsecno
  \def\newsubsec#1{\procno=0 \advance\subsecno by 1
    \edef\currentsec{\number\secno.\number\subsecno}
     \edef#1{\currentsec}}

  \newcount\appendixno \appendixno=0
  \def\newappendix#1{\procno=0 \ignorsec=0 \advance\appendixno by 1
    \ifnum\appendixno=1 \edef\appendixalpha{\hbox{A}}
      \else \ifnum\appendixno=2 \edef\appendixalpha{\hbox{B}} \fi
      \else \ifnum\appendixno=3 \edef\appendixalpha{\hbox{C}} \fi
      \else \ifnum\appendixno=4 \edef\appendixalpha{\hbox{D}} \fi
      \else \ifnum\appendixno=5 \edef\appendixalpha{\hbox{E}} \fi
      \else \ifnum\appendixno=6 \edef\appendixalpha{\hbox{F}} \fi
    \fi
    \edef#1{\appendixalpha}
    \edef\currentsec{\appendixalpha}}

  \newcount\procno \procno=0
  \def\newproc#1{\advance\procno by 1
   \ifnum\ignorsec=0 \edef#1{\currentsec.\number\procno}
                     \edef\currentproc{\currentsec.\number\procno}
   \else \edef#1{\number\procno}
         \edef\currentproc{\number\procno}
   \fi}

  \newcount\subprocno \subprocno=0
  \def\newsubproc#1{\advance\subprocno by 1
   \ifnum\subprocno=1 \edef#1{\currentproc a} \fi
   \ifnum\subprocno=2 \edef#1{\currentproc b} \fi
   \ifnum\subprocno=3 \edef#1{\currentproc c} \fi
   \ifnum\subprocno=4 \edef#1{\currentproc d} \fi
   \ifnum\subprocno=5 \edef#1{\currentproc e} \fi
   \ifnum\subprocno=6 \edef#1{\currentproc f} \fi
   \ifnum\subprocno=7 \edef#1{\currentproc g} \fi
   \ifnum\subprocno=8 \edef#1{\currentproc h} \fi
   \ifnum\subprocno=9 \edef#1{\currentproc i} \fi
   \ifnum\subprocno>9 \edef#1{TOO MANY SUBPROCS} \fi
  }

  \newcount\tagno \tagno=0
  \def\newtag#1{\advance\tagno by 1 \edef#1{\number\tagno}}



\notasec
  \newtag{\tensorexact}
  \newtag{\Afp}
  \newtag{\GpAmalgFP}
\newsec{\DVVconstr}
  \newtag{\phiax}
 \newproc{\BGNS}
 \newproc{\AContr}
  \newtag{\Visom}
  \newtag{\lambdaa}
 \newproc{\Omit}
 \newproc{\RedWords}
\newsec{\Prelims}
 \newproc{\FactorHilmod}
 \newproc{\SplitExact}
  \newtag{\threelemseq}
  \newtag{\unitalses}
 \newproc{\NoCompacts}
 \newproc{\SepSubalg}
 \newproc{\SepEmb}
  \newtag{\AfpNS}
  \newtag{\Cfp}
 \newproc{\Kexact}
 \newproc{\Folded}
  \newtag{\Bowtie}
  \newtag{\foldspanned}
 \newproc{\FoldedFacts}
 \newproc{\FoldedIdeals}
  \newtag{\JAses}
  \newtag{\IDCses}
 \newproc{\CPSplRem}
  \newtag{\stcond}
 \newproc{\capCompacts}
  \newtag{\KeuCapLeu}
\newsec{\Exactness}
 \newproc{\cpuinv}
  \newtag{\Vpkgoes}
  \newtag{\ThetaPhi}
 \newproc{\exactfp}
  \newtag{\KBAses}
  \newtag{\Calkiniso}
  \newtag{\PhikaTop}
  \newtag{\PhikaBot}
  \newtag{\Phikb}
  \newtag{\KbtA}
  \newtag{\DpDp}
  \newtag{\IIses}
  \newtag{\IoIiso}
  \newtag{\VDpVseq}
 \newproc{\ExactGps}

\newbib{\BlanchardDykemaZZEmb}
\newbib{\HarpeRobertsonValetteZZExact}
\newbib{\KirchbergZZSubalg}
\newbib{\KirchbergZZComUHF}
\newbib{\KirchbergWassermannSZZExGp}
\newbib{\LanceZZHilbertCS}
\newbib{\PisierZZExOpSp}
\newbib{\VoiculescuZZSymmetries}
\newbib{\VDNbook}
\newbib{\WassermannSZZNonExact}
\newbib{\WassermannSZZNuclEmb}
\newbib{\WassermannSZZSNU}

\topmatter
  \title 
    Exactness of reduced amalgamated free product C$^*$--algebras
  \endtitle

  \author Kenneth J\. Dykema
    \thanks Partially supported as an invited researcher funded by CNRS of
            France.
    \endthanks
  \endauthor

  \date 1 November, 1999 \enddate

  \rightheadtext{Exactness}
  \leftheadtext{Exactness}

  \address Dept.~of Mathematics,
           Texas A\&M University,
           College Station TX 77843-3368, USA
  \endaddress

  \email Ken.Dykema\@math.tamu.edu, {\it Internet URL:}
         http://www.math.tamu.edu/\~{\hskip0.1em}Ken.Dykema/
  \endemail

  \abstract
    Some completely positive maps on reduced amalgamated free products of
    C$^*$--algebras are constructed, showing that every reduced amalgamated
    free product of exact C$^*$--algebras is exact.
    Consequently, every amalgamated free product of exact discrete
    groups is exact.
  \endabstract


\endtopmatter

\document \TagsOnRight \baselineskip=18pt

\heading Introduction.\endheading
\vskip3ex

A C$^*$--algebra $A$ is said to be exact if for every short exact sequence
$$ 0\to J\to B\to B/J\to0 $$
of C$^*$--algebras and $*$--homomorphisms, the sequence
$$ 0\to A\otimes J\to A\otimes B\to A\otimes(B/J)\to0 \tag{\tensorexact} $$
of spatial tensor products is exact.
The issue was first raised when S\.~Wassermann~\cite{\WassermannSZZNonExact}
proved in an example that the sequence~(\tensorexact) need not be exact.
Later, he showed~\cite{\WassermannSZZNuclEmb} that a sufficient condition for
$A$ to be exact is that it have a nuclear embedding
$A\hookrightarrow A'$ for some C$^*$--algebra $A'$, i.e\. that for
arbitrary $\epsilon>0$ and for every finite subset $\omega\subseteq A$ there be
an integer $n$ and completely positive contractions $\Phi:A\to M_n(\Cpx)$ and
$\Psi:M_n(\Cpx)\to A'$ such that $\nm{x-\Psi\circ\Phi(x)}<\epsilon$ for
every $x\in\omega$.
The class of exact C$^*$--algebras is known to be closed under taking
subalgebras, taking spatial tensor products and taking
inductive limits.
In some remarkable papers,~\cite{\KirchbergZZSubalg}
and~\cite{\KirchbergZZComUHF}, E\.~Kirchberg proved a number of conditions
equivalent to exactness for separable, unital C$^*$--algebras, among them
nuclear embeddability.
However, in this article we will not need to make use of Kirchberg's powerful
results.
A good reference for exact C$^*$--algebras is Wassermann's
monograph~\cite{\WassermannSZZSNU}. 

In~\cite{\VoiculescuZZSymmetries}, Voiculescu introduced the noncommutative
probabilistic theory of freeness, which has turned out to be instrumental to
the study of C$^*$--algebras and von Neumann algebras associated to free
products of groups, (see for example the book~\cite{\VDNbook}).
His amalgamated or ``$B$--valued '' version of freeness is as follows:
if $A$ is a unital C$^*$--algebra with a unital C$^*$--subalgebra $B$ and a
conditional expectation (i.e\. a projection of norm $1$), $\phi:A\to B$ and if
$B\subseteq A_\iota\subseteq A$ are intermediate C$^*$--subalgebras,
($\iota\in I$), then the family $(A_\iota)_{\iota\in I}$ is said to be free
with respect to $\phi$ if $\phi(a_1a_2\cdots a_n)=0$ whenever
$a_j\in A_{\iota_j}\cap\ker\phi$ and
$\iota_1\ne\iota_2,\,\iota_2\ne\iota_3,\cdots,\iota_{n-1}\ne\iota_n$.
Voiculescu also introduced the reduced amalgamated free product of
C$^*$--algebras, which we now describe.
If $B$ is a unital C$^*$--algebra and if for a set $I$ and every $\iota\in I$,
$A_\iota$ is a unital C$^*$--algebra containing a copy of $B$ as a unital
C$^*$--subalgebra and having a conditional expectation
$\phi_\iota:A_\iota\to B$ whose GNS representation is faithful (see~\BGNS{}
below), then there is a unique unital C$^*$--algebra $A$ containing a unital
copy of $B$, with a conditional expectation $\phi:A\to B$ and with embeddings
$A_\iota\hookrightarrow A$ restricting to the identity on $B$ such that
\roster
\item"(i)" $\forall\iota\in I$ $\phi\restrict_{A_\iota}=\phi_\iota$;
\item"(ii)" the family $(A_\iota)_{\iota\in I}$ is free with respect to $\phi$;
\item"(iii)" $A$ is generated by $\bigcup_{\iota\in I}A_\iota$;
\item"(iv)" the GNS representation of $\phi$ is faithful on $A$.
\endroster
This is called the reduced amalgamated free product of C$^*$--algebras, and is
denoted by
$$ (A,\phi)=\freeprodi(A_\iota,\phi_\iota). \tag{\Afp} $$
In the case when $B=\Cpx$, the conditional expectations are just states and
the construction~(\Afp) is often called simply the reduced free product.

The following prominent example relates the reduced amalgamated free product of
C$^*$--algebras to the amalgamated free product of groups.
Let $B$ be the reduced C$^*$--algebra $C^*_\red(H)$ of a discrete group $H$,
let $G_\iota$ be a discrete group containing a copy of $H$ as a subgroup and
Take the conditional expectation $\tau_H^{G_\iota}:A_\iota\to B$ given by
$$ \tau^{G_\iota}_H(\lambda_g)=\cases
\lambda_g&\text{if }g\in H \\
0&\text{if }g\notin H. \endcases $$
Then the reduced amalgamated free product construction yields
$$ (C^*_\red(G),\tau_H^G)=\freeprodi(C^*_\red(G_\iota),\tau_H^{G_\iota})
\tag{\GpAmalgFP} $$
where $G$ is the amalgamated free product of groups $G=\fpiamalg HG_\iota$.

In this paper it is proved that in the reduced amalgamated free product of
C$^*$--algebras~(\Afp), if each $A_\iota$ is exact then $A$ is exact.
The proof proceeds by the construction of completely positive maps showing that
$A$ satisfies a condition analogous to nuclear embeddability; (see
Lemma~\FactorHilmod).
Using different techniques, Kirchberg has proved~\cite{\KirchbergZZComUHF} that
every reduced amalgamated free product of finite dimensional C$^*$--algebras is
exact.

In~\cite{\KirchbergWassermannSZZExGp}, Kirchberg and Wassermann defined a
locally compact group $G$ to be exact if the sequence
$$ 0\to J\rtimes_{\alpha\restrict_J,r}G\to A\rtimes_{\alpha,r}G\to
(A/J)\rtimes_{\tilde{\alpha},r}G\to0 $$
is exact for every continuous action $\alpha$ of $G$ on a C$^*$--algebra $A$
and every $\alpha$--invariant ideal $J$ of $A$.
They also showed that a discrete group $G$ is exact if and only if its reduced
group C$^*$--algebra $C^*_\red(G)$ is exact.
A consequence of our main result is that the class of exact discrete groups is
closed under taking amalgamated free products.

In~\S\DVVconstr, Voiculescu's construction of the reduced amalgamated free
product of C$^*$--algebras is described in detail.
In~\S\Prelims, some preliminary lemmas about exact C$^*$--algebras and  Hilbert
C$^*$--modules are proved.
In~\S\Exactness, completely positive maps on reduced amalgamated free product
C$^*$--algebras are constructed and the main result is proved.

\demo{Acknowledgements}
I would like to thank Etienne Blanchard and Marius Junge for several helpful
conversations.
Most of this work was done while I was visiting {\it l'Institut de
Math\'emat\-iques de Luminy} near Marseille.
I would like to thank the members of the department and especially \'Etienne
Blanchard and J\'er\^ome Chabert for their kind hospitality.
\enddemo

\vskip3ex
\heading\S\DVVconstr.  The construction of reduced amalgamated free products.
\endheading
\vskip3ex

In this section, by way of introducing some notation and a few conventions, we
recall Voiculescu's construction~\cite{\VoiculescuZZSymmetries} of reduced
amalgamated free products of C$^*$--algebras.

Let $B$ be a unital C$^*$--algebra,
let $I$ be a set having at least two elements and for every $\iota\in I$ let
$A_\iota$ be a unital C$^*$--algebra containing a copy of $B$ as a unital
C$^*$--subalgebra;
suppose that $\phi_\iota:A_\iota\to B$ is a conditional expectation
satisfying the property that
$$ \forall a\in A_\iota\backslash\{0\}\quad
\exists x\in A_\iota\quad\phi_\iota(x^*a^*ax)\ne0. \tag{\phiax} $$
Then
$$ (A,\phi)=\freeprodi(A_\iota,\phi_\iota) $$
denotes the reduced amalgamated free product of C$^*$--algebras, whose
construction is given below.

The construction, if $B\ne\Cpx$, depends on the theory of Hilbert
C$^*$--modules;
see Lance's book~\cite{\LanceZZHilbertCS} for a good general reference.
However, if $B=\Cpx$ then all the Hilbert C$^*$--modules considered below
become simply Hilbert spaces.

\proclaim{\BGNS}\rm
Let $E_\iota=L^2(A_\iota,\Phi_\iota)$ be the (right) Hilbert $B$--module
obtained from $A$ by separation and completion with respect to the norm
$\nm{a}=\nm{\langle a,a\rangle_{E_\iota}}^{1/2}$, where
$\langle\cdot,\cdot\rangle_{E_\iota}$ is the $B$--valued inner product,
$\langle a_1,a_2\rangle_{E_\iota}=\phi_\iota(a_1^*a_2)$;
note that this inner product is conjugate linear in the first variable, as will
be all the inner products in this paper.
We denote the map $A_\iota\to E_\iota$ arising from the definition by
$a\mapsto\widehat a$.
Let $\pi_\iota:A_\iota\to\Leu(E_\iota)$ denote the $*$--representation defined
by $\pi_\iota(a)\widehat b=\widehat{ab}$,
where as usual, for a Hilbert $B$--module $E$, $\Leu(E)$ denotes the
C$^*$--algebra of all adjointable bounded $B$--module operators on $E$; we will
also let $\Keu(E)$ denote the C$^*$--subalgebra (in fact, the ideal) of
$\Leu(E)$ generated by the collection of all operators of the form
$\theta_{x,y}$ for $x,y\in E$, given by
$\theta_{x,y}(e)=x\langle y,e\rangle_E$.
Note that the condition~(\phiax) is equivalent to faithfulness of
$\pi_\iota$.
Consider the specified element $\xi_\iota=\widehat{1_{A_\iota}}\in E_\iota$.
We will call $(\pi_\iota,E_\iota,\xi_\iota)$ the GNS representation of
$(A_\iota,\phi_\iota)$ and write
$(\pi_\iota,E_\iota,\xi_\iota)=\GNS(A_\iota,\phi_\iota)$, though this is
technically a misnomer unless $B=\Cpx$.
\endproclaim

\proclaim{\AContr}\rm
Voiculescu's construction of $A$ proceeds via the construction of a Hilbert
$B$--module $E$, which can later be seen to be $L^2(A,\phi)$.
Note that the subspace $\xi_\iota B$ is a complemented submodule of $E_\iota$
that is invariant under the left action $\pi_\iota\restrict_B$ of $B$;
indeed, $\theta_{\xi_\iota,\xi_\iota}$ is the projection onto $\xi_\iota B$.
We will denote the complementing submodule by
$E_\iota\oup=P_\iota\oup E_\iota$,
where we define $P_\iota\oup=1-\theta_{\xi_\iota,\xi_\iota}\in\Leu(E_\iota)$.
Let
$$ E=\xi B\oplus\bigoplus
\Sb n\in\Nats \\ \iota_1,\ldots,\iota_n\in I\\
\iota_1\ne\iota_2,\,\iota_2\ne\iota_3,\,\ldots,\,\iota_{n-1}\ne\iota_n\endSb
E_{\iota_1}\oup\otimes_BE_{\iota_2}\oup\otdts BE_{\iota_n}\oup; $$
here $\xi B$ denotes the C$^*$--algebra $B$ considered as a Hilbert $B$--module
and with specified element $\xi=1_B$,
the tensor products are internal tensor products arising in relation to the
$*$--homomorphisms $P_\iota\oup\pi_\iota\restrict_B(\cdot)P_\iota\oup$ from $B$
to $\Leu(E_\iota\oup)$,
and $\Nats$ for us always  means the positive integers excluding $0$.
All of the tensor products denoted $\otimes_B$ in this
paper will be internal tensor products defined with respect to some
$*$--homomorphisms from $B$ into $\Leu(F)$ for the various Hilbert $B$--modules
$F$, these $*$--homomorphisms arising canonically from $\pi\restrict_B$ and the
several constructions employed on Hilbert $B$--modules.
The Hilbert $B$--module $E$ constructed above is call the {\it free product}
of the $E_\iota$ with respect to specified vectors $\xi_\iota$, and will be
denoted by $(E,\xi)=\freeprodi(E_\iota,\xi_\iota)$.

For $\iota\in I$ let
$$ E(\iota)=\eta_\iota B\oplus\bigoplus
\Sb n\in\Nats \\ \iota_1,\ldots,\iota_n\in I\\
 \iota_1\ne\iota_2,\,\ldots,\,
 \iota_{n-1}\ne\iota_n \\
 \iota_1\ne\iota \endSb
E_{\iota_1}\oup\otimes_BE_{\iota_2}\oup\otdts BE_{\iota_n}\oup,
$$
where $\eta_\iota B$ is a copy of the Hilbert $B$--module $B$ with
$\eta_\iota=1_B$, and let
$$ V_\iota:E_\iota\otimes_BE(\iota)\to E \tag{\Visom} $$
be the unitary operator defined as follows.
In order to distinguish the tensor product in~(\Visom) from those
appearing in elements of $E$ and $E(\iota)$, we will use the symbol
$\otimesdd$ for it;
then $V_\iota$ is given by
$$ \align
V_\iota:\,&\xi_\iota\otimesdd\eta_\iota\mapsto\xi \\
&\zeta\otimesdd\eta_\iota\mapsto\zeta \\
&\xi_\iota\otimesdd(\zeta_1\otdt\zeta_n)\mapsto\zeta_1\otdt\zeta_n \\
&\zeta\otimesdd(\zeta_1\otdt\zeta_n)
 \mapsto\zeta\otimes\zeta_1\otdt\zeta_n,
\endalign $$
whenever $\zeta\in E_\iota\oup$ and $\zeta_j\in E_{\iota_j}\oup$ with
$\iota\ne\iota_1,\,\iota_1\ne\iota_2,\,\ldots,\,\iota_{n-1}\ne\iota_n$.
Let $\lambda_\iota:A_\iota\to\Leu(E)$ be the $*$--homomorphism given by
$$ \lambda_\iota(a)=V_\iota(\pi_\iota(a)\otimes1)V_\iota^*. $$
Then $A$ is defined to be the C$^*$--algebra generated by
$\bigcup_{\iota\in I}\lambda_\iota(A_\iota)$, and $\phi:A\to B$ is the
conditional expectation $\phi(\cdot)=\langle\xi,\,\cdot\,\xi\rangle_E$.
It is important to note that for $b\in B$, the operator $\lambda_\iota(b)$ on
$E$ does not depend on $\iota$.

Let $A_\iota\oup=A_\iota\cap\ker\phi_\iota$.
Note that if
$a\in A_\iota\oup$
and if $\zeta_j\in E_{\iota_j}\oup$ for $\iota_1,\ldots,\iota_n\in I$, $n\ge2$,
and $\iota_j\ne\iota_{j+1}$, then
$$ \lambda_\iota(a)(\zeta_1\otdt\zeta_n)=\cases
\widehat a\otimes\zeta_1\otimes\zeta_2\otdt\zeta_n&\text{if }\iota\ne\iota_1 \\
 \vspace{2ex}
\aligned (a\zeta_1-\xi_{\iota_1}\langle\xi_{\iota_1},a\zeta_1\rangle)
 \otimes&\zeta_2\otdt\zeta_n \\
 +\;\pi_{\iota_2}(\langle\xi_{\iota_1},a\zeta_1\rangle)
 &\zeta_2\otdt\zeta_n\endaligned
 &\text{if }\iota=\iota_1.
\endcases \tag{\lambdaa} $$
\endproclaim

\proclaim{\Omit}\rm
Throughout this paper we will usually omit to write the representations
$\lambda_\iota$ and $\pi_\iota$, and we will simply
think of the $A_\iota$ as C$^*$--subalgebras of $\Leu(E_\iota)$ when
it suits us and of $\Leu(E)$ when it suits us,
and of $B$ acting on the left and the right of just about everything.
\endproclaim

\proclaim{\RedWords}\rm
The free product C$^*$--algebra
$A$ is the closed linear span of $B$ together with the set of
all {\it reduced words} of the form $w=a_1a_2\cdots a_q$ where $q\in\Nats$,
$a_j\in A_{\iota_j}\oup$, $\iota_1,\ldots,\iota_q\in I$ and
$\iota_1\ne\iota_2,\,\ldots,\iota_{q-1}\ne\iota_q$.
The {\it length} of this word $w$ is $q$.
Using~(\lambdaa) one can see that if $w=a_1\cdots a_q$ is a reduced word of
length $q$ and if
$\zeta_1\otdt\zeta_n\in E_{\iota'_1}\oup\otdts BE_{\iota'_n}\oup\subseteq E$
for $n>q$, then when acting on $\zeta_1\otdt\zeta_n$, $w$ only sees the first
bit, $\zeta_1\otdt\zeta_q$, and $w(\zeta_1\otdt\zeta_n)$ is a linear
combination of simple tensors having lengths between $n-q$ and $n+q$, and each
with the same tail $\cdots\otimes\zeta_{q+1}\otdt\zeta_n$, (possibly multiplied
on the left by an element of $B$).
The tail simply hangs on for the ride, so to speak.
\endproclaim

\vskip3ex
\heading\S\Prelims.  Preliminary lemmas and a construction \endheading
\vskip3ex

In this section are assembled a few preliminary results and a construction
involving Hilbert C$^*$--modules.
We begin with an easy generalization of the proof of Wassermann's
result~\cite{\WassermannSZZNuclEmb} that nuclear embeddability implies
exactness.
\proclaim{Lemma \FactorHilmod}
Let $A'$ be a C$^*$--algebra and let $A$ be a C$^*$--subalgebra of
$A'$.
Suppose that for every finite subset $\omega\subseteq A$ and every
$\epsilon>0$ there is an exact C$^*$--algebra $D$ and there are
completely positive contractions, $\Phi:A\to D$ and $\Psi:D\to A'$, such that
$\nm{x-\Psi\circ\Phi(x)}<\epsilon$ for every $x\in\omega$.
Then $A$ is exact.
\endproclaim
\demo{Proof}
Using the hypotheses and taking the directed set of all finite subsets of $A$,
one constructs a net $(\Phi_\lambda,D_\lambda,\Psi_\lambda)$ of exact
C$^*$--algebras $D_\lambda$ and completely positive maps
$\Phi_\lambda:A\to D_\lambda$ and $\Psi_\lambda:D_\lambda\to A'$ such that
$$ \forall x\in A\qquad
\lim_\lambda\nm{\Psi_\lambda\circ\Phi_\lambda(x)-x}=0. $$
Let
$$ 0\to I\to B\overset\pi\to\to C\to0 $$
be an exact sequence of C$^*$--algebras and $*$--homomorphisms.
Let $x\in(A\otimes B)\cap\ker(\id_A\otimes\pi)$.
We will show that $x\in A\otimes I$, which will prove the lemma.
For every $\lambda$ we have
$$ 0=(\Phi_\lambda\otimes\id_C)\circ(\id_A\otimes\pi)(x)
=(\id_{D_\lambda}\otimes\pi)\circ(\Phi_\lambda\otimes\id_B)(x), $$
so
$(\Phi_\lambda\otimes\id_B)(x)\in(D_\lambda\otimes B)
\cap\ker(\id_{D_\lambda}\otimes\pi)$.
Since $D_\lambda$ is exact, this implies that
$(\Phi_\lambda\otimes\id_B)(x)\in D_\lambda\otimes I$.
Hence
 $\bigl((\Psi_\lambda\circ\Phi_\lambda)\otimes\id_B\bigr)(x)\in A'\otimes I$.
We then have that
$$ x=\lim_\lambda\bigl((\Psi_\lambda\circ\Phi_\lambda)\otimes\id_B\bigr)(x)
\in A'\otimes I. $$
But it is easily seen, using an approximate identity for $I$, that
$(A\otimes B)\cap(A'\otimes I)=A\otimes I$.
Therefore we find that $x\in A\otimes I$, as required.
\QED

The next result is well known to experts;
proofs can be found at~\scite{\HarpeRobertsonValetteZZExact}{Prop\. 2}
and~\scite{\KirchbergZZComUHF}{7.1}, and another proof is possible using exact
operator spaces~\cite{\PisierZZExOpSp}.
I would like to thank Etienne Blanchard for first bringing it to my attention
and showing me a proof.
Here and everywhere in this paper, ideals of C$^*$--algebras are closed,
two--sided ideals.

\proclaim{Lemma \SplitExact}
Let $A$ be a C$^*$--algebra, let $J$ be an ideal of $A$ and let $q:A\to A/J$ be
the quotient map.
Suppose that the short exact sequence
$$ 0\to J\to A\overset q\to\to A/J\to 0 $$
has a completely positive contractive splitting, i.e\. a completely positive
contraction, $s:A/J\to A$, such that $q\circ s=\id_{A/J}$.
If $A/J$ and $J$ are exact C$^*$--algebras, then $A$ is an exact
C$^*$--algebra.
\endproclaim

\proclaim{Lemma \NoCompacts}
Let $A$ be a unital C$^*$--algebra having a unital C$^*$--subalgebra $B$ and a
conditional expectation $\phi:A\to B$.
Let $C(\Tcirc)$ be the C$^*$--algebra of all continuous functions on the circle
and let $\tau$ be the tracial state on $C(\Tcirc)$ given by integration with
respect to Haar measure.
Consider the conditional expectation
$$ \phi\otimes\tau:A\otimes C(\Tcirc)\to B\otimes1\cong B $$
and let $(\pi,F,\xi)$ be the GNS representation of
$(A\otimes C(\Tcirc),\phi\otimes\tau)$ as described in~\BGNS.
Then
$$ \pi\bigl(A\otimes C(\Tcirc)\bigr)\cap\Keu(F)=\{0\}. $$
\endproclaim
\demo{Proof}
Let
$$ \align
(\sigma,E,\eta)&=\GNS(A,\phi) \\
(\rho,\Hil,\eta')&=\GNS(C(\Tcirc),\tau);
\endalign $$
then $\Hil=L^2(\Tcirc)$ and $\rho(f)$, for $f\in C(\Tcirc)$, is multiplication
by $f$.
The Hilbert $A\otimes C(\Tcirc)$--module $F$ is canonically isomorphic to the
external tensor product $E\otimes\Hil$, $\pi$ is thereby identified with
$\sigma\otimes\rho$ and $\Keu(F)$ is identified with
$\Keu(E)\otimes\Keu(\Hil)$.
But $\pi(A\otimes C(\Tcirc))$ commutes with $1\otimes\rho(C(\Tcirc))$, and an
elementary argument shows that
$$ \bigl(\Leu(E)\otimes\Keu(\Hil)\bigr)
\cap\bigl(1\otimes\rho(C(\Tcirc))\bigr)'=\{0\}. $$
\QED

\proclaim{Lemma~\SepSubalg}
Let $B$ be a C$^*$--algebra and let $I$ be a countable set.
For every $\iota\in I$ let $A_\iota$ be a C$^*$--algebra containing $B$ as a
C$^*$--subalgebra and having a conditional expectation
$\phi_\iota:A_\iota\to B$.
Let $Y$ be a countable subset of $B$ and for every $\iota\in I$ let $X_\iota$
be a countable subset of $A_\iota$.
Then there are separable C$^*$--subalgebras $D\subseteq B$ and
$C_\iota\subseteq A_\iota$ ($\iota\in I$), such that $Y\subseteq D$,
$X_\iota\subseteq C_\iota$, $D\subseteq C_\iota$ and
$\phi_\iota(C_\iota)=D$ for every $\iota\in I$.
\endproclaim
\demo{Proof}
Let $D_1=C^*(Y)\subseteq B$ and $C_{\iota,1}=C^*(X_\iota)\subseteq A_\iota$, and
for $n\in\Nats$ define recursively
$$ C_{\iota,n+1}=C^*(C_{\iota,n}\cup D_n)\qquad\text{and}\qquad
D_{n+1}
=C^*\biggl(D_n\cup\bigcup_{\iota\in I}\phi_\iota(C_{\iota,n+1})\biggr). $$
Finally, let
$$ C_\iota=\overline{\bigcup_{n\ge1}C_{\iota_n}}
\qquad\text{and}\qquad D=\overline{\bigcup_{n\ge1}D_n}. $$
\QED

\proclaim{Remark \SepEmb}\rm
If in Lemma~\SepSubalg each $A_\iota$ is unital, having $B$ as a unital
subalgebra and if the GNS representation of each $\phi_\iota$ is faithful, then
consider the reduced amalgamated free product
$$ (A,\phi)=\freeprodi(A_\iota,\phi_\iota). \tag{\AfpNS} $$
For $D$ and $C_\iota$ as constructed in the lemma, take the conditional
expectations $\psi_\iota\eqdef\phi_\iota\restrict_{C_\iota}:C_\iota\to D$ and
the reduced amalgamated free product
$$ (C,\phi)=\freeprodi(C_\iota,\psi_\iota). \tag{\Cfp} $$
Then by the main result of~\cite{\BlanchardDykemaZZEmb}, the embeddings
$C_\iota\hookrightarrow A_\iota$ extend to an embedding $C\hookrightarrow A$.
From this it is easily seen that every reduced amalgamated free product
C$^*$--algebra $A$ as in~(\AfpNS) is the inductive limit of reduced amalgamated
free products of separable subalgebras of the $A_\iota$ as in~(\Cfp).
\endproclaim

\proclaim{Lemma \Kexact}
Let $B$ be an exact C$^*$--algebra and let $E$ be a countably generated Hilbert
$B$--module.
Then $\Keu(E)$ is an exact C$^*$--algebra.
\endproclaim
\demo{Proof}
By Kasparov's stabilization lemma, $E$ is a complemented C$^*$--submodule of
$\Hil_B\eqdef\Hil\otimes B$, the external tensor product of a separable
infinite dimensional Hilbert space $\Hil$ and the Hilbert $B$--module $B$.
Thus $\Keu(E)$ is a C$^*$--subalgebra of
$\Keu(\Hil_B)\cong\Keu(\Hil)\otimes B$, which is exact.
\QED

For the rest of this section we concentrate on the following construction.
\proclaim{Definition \Folded}\rm
Let $B_1$ and $B_2$ be C$^*$--algebras and let $E_i$ be a Hilbert
$B_i$--module, ($i=1,2$);
let $\pi:B_1\to\Leu(E_2)$ be a $*$--homomorphism and let $E_1\otimes_\pi E_2$
be the associated internal tensor product.
If $A\subseteq\Leu(E_2)$ is a C$^*$--subalgebra such that
$\pi(B_1)A\subseteq A$, then let
$$ \Keu(E_1)\bowtie_\pi A \tag{\Bowtie} $$
be the C$^*$--subalgebra of $\Leu(E_1\otimes_\pi E_2)$ generated by
$$ \{\theta_ea\theta_{e'}^*\mid e,e'\in E_1,\,a\in A\}, \tag{\foldspanned} $$
where $\theta_e\in\Leu(E_2,E_1\otimes_\pi E_2)$ is the operator defined by
$\theta_e(f)=e\otimes f$, (see~\scite{\LanceZZHilbertCS}{4.6}).
\endproclaim

By identifying $E_1\otimes_\pi E_2$ with $E_1\otimes_\pi A\otimes_AE_2$, where
we regard $\pi$ as also a $*$--homomorphism from $B_1$ into $\Leu(A)$ for the
Hilbert $A$--module $A$, one sees that $\Keu(E_1)\bowtie_\pi A$ is canonically
isomorphic to $\Keu(E_1\otimes_\pi A)\otimes1_{E_2}$.
However, we persist with the notation~(\Bowtie) because it seems more
convenient for keeping track of the tensor product structure of the Hilbert
modules.
It may be helpful to realize that if $B_1=\Cpx$ and $\pi$ is unital then
$\Keu(E_1)\bowtie_\pi A$ is simply $\Keu(E_1)\otimes A$.

Some easy facts about this construction are collected in the lemma below.
\proclaim{Lemma \FoldedFacts}
\roster
\item"(i)" $\Keu(E_1)\bowtie_\pi A$ is the closed linear span of the
set~(\foldspanned).
\item"(ii)" $\Keu(E_1)\bowtie_\pi\Keu(E_2)=\Keu(E_1\otimes_\pi E_2)$.
\item"(iii)" Let $\Hil_{B_1}$ be the Hilbert $B_1$--module that is the external
tensor product $\Hil\otimes B_1$, where $\Hil$ is a separable, infinite
dimensional Hilbert space.
Then $\Keu(\Hil_{B_1})\bowtie_\pi A$ is canonically isomorphic to
$\Keu(\Hil)\otimes\pi(B_1)A\pi(B_1)$.
\item"(iv)" Suppose $E_1$ is countably generated.
Then by the Kasparov stabilization theorem, $E_1$ is isomorphic to a
complemented submodule of $\Hil_{B_1}$.
Let $P$ be the projection from $\Hil_{B_1}$ onto $E_1$.
The inclusion $E_1\hookrightarrow\Hil_{B_1}$ provides an inclusion
$\Keu(E_1)\bowtie_\pi A\hookrightarrow\Keu(\Hil_{B_1})\bowtie_\pi A$ and the
map $\Phi(x)=(P\otimes1)x(P\otimes1)$ is a conditional expectation from
$\Keu(\Hil_{B_1})\bowtie_\pi A$ onto $\Keu(E_1)\bowtie_\pi A$.
\endroster
\endproclaim
\demo{Proof}
(i) follows from $\theta_{e'}^*\theta_e=\pi(\langle e',e\rangle)$.

(ii) follows from
$\theta_e\theta_{f,f'}\theta_{e'}^*=\theta_{e\otimes f,e'\otimes f'}$.

In~(iii), the isomorphism is
$\rho:\Keu(\Hil_{B_1})\bowtie_\pi A\to\Keu(\Hil)\otimes\pi(B_1)A\pi(B_1)$ given
by, for $v,v'\in\Hil$, $b,b'\in B$ and $a\in A$,
$$ \rho(\theta_{v\otimes b}a\theta_{v'\otimes b'})
=\theta_{v,v'}\otimes\pi(b)a\pi(b'). $$
\QED

Perhaps slightly less obvious facts about this construction are below.
\proclaim{Lemma \FoldedIdeals}
Let $J$ be an ideal of $A$ and suppose that $E_1$ is countably generated;
note that $\pi_1(B_1)J\subseteq J$.
\roster
\item"(i)" Then $\Keu(E_1)\bowtie_\pi J$ is an ideal of
$\Keu(E_1)\bowtie_\pi A$ and the quotient C$^*$--algebra
$$ \frac{\Keu(E_1)\bowtie_\pi A}{\Keu(E_1)\bowtie_\pi J} $$
is isomorphic to a C$^*$--subalgebra of $\Keu\otimes(A/J)$, where $\Keu$ is the
algebra of all compact operators on a separable infinite dimensional Hilbert
space.
\item"(ii)" If the short exact sequence
$$ 0\to J\to A\to A/J\to 0 \tag{\JAses} $$
has a completely positive contractive splitting $s:(A/J)\to A$ then the short
exact sequence
$$ 0\to\Keu(E_1)\bowtie_\pi J\to\Keu(E_1)\bowtie_\pi A\to
\frac{\Keu(E_1)\bowtie_\pi A}{\Keu(E_1)\bowtie_\pi J}\to0 $$
has a completely positive contractive splitting
$$ \st:\frac{\Keu(E_1)\bowtie_\pi A}{\Keu(E_1)\bowtie_\pi J}\to
\Keu(E_1)\bowtie_\pi A. $$
\endroster
\endproclaim
\demo{Proof}
Let us write $D=\Keu(E_1)\bowtie_\pi A$ and $I=\Keu(E_1)\bowtie_\pi J$.
It is clear that $I$ is an ideal of $D$.
Let $\rho:\Keu(\Hil_{B_1})\bowtie_\pi A\to\Keu\otimes A$ be the isomorphism
indicated in Lemma~\FoldedFacts(iii) and let
$\Phi:\Keu(\Hil_{B_1})\bowtie_\pi A\to\Keu(E_1)\bowtie_\pi A$ be the
conditional expectation in~\FoldedFacts(iv).
One easily verifies that
$\rho(\Keu(\Hil_{B_1})\bowtie_\pi J)=\Keu\otimes J\subseteq\Keu\otimes A$
and $\Phi(\KHil(\Hil_{B_1})\bowtie_\pi J)=I$.
Therefore we see that $\rho(D)\cap(\Keu\otimes J)=\rho(I)$, and hence that
$D/I$ is isomorphic to the image, call it $C$, of $\rho(D)$ in the quotient
$(\Keu\otimes A)/(\Keu\otimes J)\cong\Keu\otimes(A/J)$.
This proves~(i).

If $s:A/J\to A$ is a completely positive contractive lifting of the short exact
sequence~(\JAses), then $\id_\Keu\otimes s:\Keu\otimes(A/J)\to\Keu\otimes A$ is
a completely positive contraction.
Let $q_D:\rho(D)\to C$ be the quotient map;
we seek a completely positive splitting for the short exact sequence
$$ 0\to I\to D\overset{q_D\circ\rho}\to\longrightarrow C\to0. \tag{\IDCses} $$
Let $x\in C\subseteq\Keu\otimes(A/J)$ and let $y=(\id_\Keu\otimes s)(x)$.
Then $y\in\rho(D)+(\Keu\otimes J)$;
hence $\rho\circ\Phi\circ\rho^{-1}(y)-y\in\Keu\otimes J$.
But $\rho\circ\Phi\circ\rho^{-1}(y)\in\rho(D)$ and
$q_D\circ\rho\circ\Phi\circ\rho^{-1}(y)=x$.
Hence
$$ \st=\Phi\circ\rho^{-1}\circ(\id_\Keu\otimes s)\restrict_C:C\to D $$
is the desired completely positive contractive splitting of~(\IDCses).
\QED

\proclaim{Remark~\CPSplRem}\rm
It is natural to ask whether the splitting $\st$ constructed above satisfies
$$ \st\bigl(\theta_ea\theta_{e'}^*+\Keu(E_1\bowtie_\pi J)\bigr)
=\theta_es(a+J)\theta_{e'} \tag{\stcond} $$
for every $a\in A$ and $e,e'\in E_1$.
I don't know the answer to this question in general, but if $s$ is assumed to
satisfy the additional condition
$s(\pi(b)a+J)=\pi(b)s(a+J)$ then it is straightforward to show that~(\stcond)
holds.
\endproclaim

\proclaim{Lemma \capCompacts}
We have
$$ \bigl(\Keu(E_1)\bowtie_\pi\Leu(E_2)\bigr)\bigcap
\bigl(\Leu(E_1)\otimes1_{E_2}\bigr)=\Keu(E_1)\otimes1_{E_2}.
\tag{\KeuCapLeu} $$
\endproclaim
\demo{Proof}
Let $(u_\lambda)_{\lambda\in\Lambda}$ be an approximate identity for
$\Keu(E_1)$, where each $u_\lambda$ is
of the form $\sum_{i=1}^n\theta_{e_i,e_i'}$, ($n\in\Nats$, $e_i,e_i'\in E_1$).
Note that $\theta_{e_i}1_{E_2}\theta_{e_i'}^*=\theta_{e_i,e_i'}\otimes1_{E_2}$,
and $(u_\lambda\otimes1_{E_2})_{\lambda\in\Lambda}$ is an approximate identity
for $\Keu(E_1)\bowtie_\pi\Leu(E_2)$.
If $y$ belongs to the left--hand--side of~(\KeuCapLeu) then $y=x\otimes1_{E_2}$
for some $x\in\Leu(E_1)$.
But also
$y=\lim_{\lambda}(u_\lambda xu_\lambda\otimes1_{E_2})
\in\Keu(E_1)\otimes1_{E_2}$.
\QED

\vskip3ex
\heading\S\Exactness.  Exactness of free product C$^*$--algebras \endheading
\vskip3ex

This section contains the proof of the main theorem, that every reduced
amalgamated free product of exact C$^*$--algebras is exact.
Let
$$ (A,\phi)=\freeprodi(A_\iota,\phi_\iota) $$
be a reduced amalgamated free product of C$^*$--algebras and let $E$ be the
Hilbert $B$--module as described in~\S\DVVconstr.
Consider the submodules of $E$,
$$ \Eto k=\xi B\oplus\bigoplus
\Sb n\in\{1,2,\ldots,k\} \\ \iota_1,\ldots,\iota_n\in I\\
\iota_1\ne\iota_2,\iota_2\ne\iota_3,\ldots,\iota_{n-1}\ne\iota_n\endSb
E_{\iota_1}\oup\otimes_BE_{\iota_2}\oup\otdts BE_{\iota_n}\oup,
\qquad(k\in\Nats) $$
and $\Eto0=\xi B$.
Let $\Pto k$ denote the projection from $E$ onto $\Eto k$.
Consider the completely positive unital maps
$$ \Phi_k:\Leu(E)\to\Leu(\Eto k) $$
obtained by compressing:
$\Phi_k(x)=\Pto kx\restrict_{\Eto k}$.
In the following lemma we will find completely positive unital maps going the
other way which are approximately left inverses on elements of
$A\subseteq\Leu(E)$.

\proclaim{Lemma \cpuinv}
There are completely positive unital maps, $\Psi_k:\Leu(\Eto k)\to\Leu(E)$
such that $\lim_{k\to\infty}\nm{a-\Psi_k\circ\Phi_k(a)}=0$ for every $a\in A$.
\endproclaim
\demo{Proof}
We will first define, for every two positive integers $p$ and $k$ with $p<k$, an
isometry
$$ V_{p,k}:E\to\Eto k\otimes_BE \tag{\Vpkgoes} $$
belonging to $\Leu(E,\Eto k\otimes_BE)$.
Let us use the symbol $\otimesdd$ for the tensor product
in~(\Vpkgoes), in order to distinguish it from tensor products appearing in
elements of $\Eto k$ or $E$.
We let $V_{p,k}(\xi)=\xi\otimesdd\xi$ and
$$ V_{p,k}(\zeta_1\otdt\zeta_n)=\cases
(\zeta_1\otdt\zeta_n)\otimesdd\xi&\text{if }1\le n\le p \\ \vspace{2ex}
\left(\aligned\sum_{j=0}^{n-p-1}&\frac1{\sqrt{k-p}}
 (\zeta_1\otdt\zeta_{p+j})\otimesdd(\zeta_{p+j+1}\otdt\zeta_n) \\
 +\;\;&\frac{\sqrt{k-n}}{\sqrt{k-p}}(\zeta_1\otdt\zeta_n)\otimesdd\xi
 \endaligned\right) 
 &\text{if }p<n\le k \\ \vspace{2ex}
\dsize\sum_{j=0}^{k-p-1}\frac1{\sqrt{k-p}}(\zeta_1\otdt\zeta_{p+j})
 \otimesdd(\zeta_{p+j+1}\otdt\zeta_n)
 &\text{if }k<n, \endcases $$
whenever $n\in\Nats$ and $\zeta_j\in E_{\iota_j}\oup$,
for some $\iota_1,\ldots,\iota_n\in I$ with
$\iota_1\ne\iota_2,\ldots,\iota_{n-1}\ne\iota_n$.
It is not difficult to check that $V_{p,k}\in \Leu(E,\Eto k\otimes_BE)$ and is an
isometry.

Let $\Theta_{p,k}:\Leu(\Eto k)\to\Leu(E)$ be the completely positive
unital map defined by
$$ \Theta_{p,k}(x)=V_{p,k}^*(x\otimes1)V_{p,k}. $$
We claim that
$$ \nm{a-\Theta_{p,k}\circ\Phi_k(a)}\to0 \tag{\ThetaPhi} $$
for every $a\in A$, as $k\to\infty$ and $p\to\infty$ in such a way that
$(k-p)\to\infty$.
It is straightforward to show that $\Theta_{p,k}\circ\Phi_k(b)=b$ for every
$b\in B\subseteq A$;
hence, to show~(\ThetaPhi) for every $a\in A$, it will suffice to show it for
every reduced word, $a=a_1\cdots a_q$, (see~\RedWords).

Let us now consider some more submodules of $E$;
for $n\in\Nats$ let
$$ E_{(n)}=\bigoplus\Sb \iota_1,\ldots,\iota_n\in I\\
\iota_1\ne\iota_2,\ldots,\iota_{n-1}\ne\iota_n\endSb
E_{\iota_1}\oup\otimes_BE_{\iota_2}\oup\otdts BE_{\iota_n}\oup, $$
and let $E_{(0)}=\xi B$;
let $P_{(n)}$ be the projection from $E$ onto $E_{(n)}$.
We can think of operators, $x\in\Leu(E)$, as infinite matrices indexed by
$\{0\}\cup\Nats$, where the $(n,m)$th entry is $P_{(n)}xP_{(m)}$.
We will let $S_d(x)$ be the matrix consisting of only the $d$th diagonal:
$$ S_d(x)=\sum_{n=\max(0,-d)}^\infty P_{(n+d)}xP_{(n)}, $$
where the sum converges in the strict topology.
Note that $\nm{S_d(x)}\le\nm x$.
Let $a=a_1\cdots a_q$ be a reduced word of length $q$,
$a_j\in A_{\iota'_j}\oup$;
then $S_d(a)$ vanishes whenever $|d|>q$; thus $a=\sum_{d=-q}^qS_d(a)$.
We will show that for every $d\in\{-q,-q+1,\ldots,q\}$,
$$ \nm{S_d(a)-\Theta_{p,k}\circ\Phi_k(S_d(a))}\to0 $$
as $k$, $p$, and $k-p$ all tend to infinity, which will suffice to
prove~(\ThetaPhi).
Fix $d$ and write $y=S_d(a)$ for convenience.
We may and do assume that $p>2q$ and $k-p>2q$.
Checking the several cases, one shows that $y-\Theta_{p,k}\circ\Phi_k(y)=yR$
where $R\in\Leu(E)$ is the operator that multiplies $\xi$ by $0$ and every
element in $E_{(n)}$ by a real number $R_n$, given by
$$ R_n=\cases
0&\text{if }n\le p\text{ and }n+d\le p \\
\left(1-\sqrt{\tfrac{k-n-d}{k-p}}\right)
 &\text{if }n\le p<n+d\le k \\
\left(1-\sqrt{\tfrac{k-n}{k-p}}\right)
 &\text{if }n+d\le p<n\le k \\
\left(\frac{k-n-\min(0,d)-(k-n)^{1/2}(k-n-d)^{1/2}}{k-p}\right)
 &\text{if }p<n\le k\text{ and }p<n+d\le k \\
\left(\frac d{k-p}\right)
 &\text{if }p<n\le k<n+d \\
\left(\frac{-d-\sqrt{k-n-d}}{k-p}\right)
 &\text{if }p<n+d\le k<n \\
\left(\frac{|d|}{k-p}\right)
 &\text{if }k<n\text{ and }k<n+d.
\endcases $$
Let us write down this computation in the case where $p<n\le k$ and
$p<n+d\le k$.
Let $\iota_1,\ldots,\iota_n\in I$ satisfy
$\iota_1\ne\iota_2,\,\ldots,\,\iota_{n-1}\ne\iota_n$ and let
$\zeta_j\in E_{\iota_j}\oup$.
We will investigate
$$ \bigl(y-\Theta_{p,k}\circ\Phi_k(y)\bigr)(\zeta_1\otdt\zeta_n)
=\Bigl(y-V_{p,k}^*(\Phi_k(y)\otimes1)V_{p,k}\Bigr)(\zeta_1\otdt\zeta_n). $$
We have
$$ \align
V_{p,k}(\zeta_1\otdt\zeta_n)
&=\left(\sum_{j=0}^{n-p-1}\tfrac1{\sqrt{k-p}}(\zeta_1\otdt\zeta_{p+j})
\otimesdd(\zeta_{p+j+1}\otdt\zeta_n)\right) \\ \vspace{2ex}
&\;+\tfrac{\sqrt{k-n}}{\sqrt{k-p}}(\zeta_1\otdt\zeta_n)\otimesdd\xi
\endalign $$
and for every $j\in\{0,1,\ldots,n-p\}$,
$$ \Phi_k(y)(\zeta_1\otdt\zeta_{p+j})=y(\zeta_1\otdt\zeta_{p+j})
\in E_{(p+j+d)}. $$
If $0\le j<-d$ then $p+j+d<p$, so
$$ V_{p,k}^*\bigl(y(\zeta_1\otdt\zeta_{p+j})
\otimesdd(\zeta_{p+j+1}\otdt\zeta_n)\bigr)=0. $$
If $\max(-d,0)\le j<n-p$ then since $p>q$ and using~\RedWords{} we have
$$ \align
V_{p,k}^*&\bigl(y(\zeta_1\otdt\zeta_{p+j})
 \otimesdd(\zeta_{p+j+1}\otdt\zeta_n)\bigr)= \\
&=\tfrac1{\sqrt{k-p}}\bigl(y(\zeta_1\otdt\zeta_{p+j})
 \otimes(\zeta_{p+j+1}\otdt\zeta_n)\bigr) \\
&=\tfrac1{\sqrt{k-p}}y(\zeta_1\otdt\zeta_n).
\endalign $$
Finally,
$$ V_{p,k}^*\bigl(y(\zeta_1\otdt\zeta_n)\otimesdd\xi\bigr)
=\tfrac{\sqrt{k-n-d}}{\sqrt{k-p}}y(\zeta_1\otdt\zeta_n). $$
Putting this together we get
$$ \bigl(y-\Theta_{p,k}\circ\Phi_k(y)\bigr)\restrict_{E_{(n)}}
=\left(\frac{k-n-\min(0,d)-(k-n)^{1/2}(k-n-d)^{1/2}}{k-p}\right)
y\restrict_{E_{(n)}}, $$
as required.
The other six cases are similar but a bit easier to check.

Since $|d|\le q$ and $q$ is fixed, we see that $R_n$ tends to zero uniformly in
$n$ as $p$, $k$ and $k-p$ all tend to infinity.
We have shown~(\ThetaPhi).
Letting $\Psi_k=\Theta_{[k/2],k}$ finishes the proof.
\QED

In the reduced amalgamated free product
$(A,\phi)=\freeprodi(A_\iota,\phi_\iota)$ from the above lemma, if we assume
that $I$ is finite and each of the C$^*$--algebras $A_\iota$ is finite
dimensional, then each of the submodules $\Eto k$ is finite dimensional;
it follows that $A$ is nuclearly embeddable, and hence is exact;
this allows an alternative proof of Kirchberg's
result~\scite{\KirchbergZZComUHF}{7.2} that reduced amalgamated free products
of finite dimensional C$^*$--algebras are exact.
In order to prove our main result, we will use the notion analogous to
nuclear embeddability as found in Lemma~\FactorHilmod{} and we will prove, for
all $k$, exactness of a C$^*$--algebra containing $\Phi_k(A)$.

\proclaim{Theorem \exactfp}
Suppose that $B$ is a unital exact C$^*$--algebra, $I$ is a set and
for every $\iota\in I$ $A_\iota$ is a unital exact C$^*$--algebra
containing $B$ as a unital C$^*$--subalgebra and having a conditional
expectation, $\phi_\iota$, from $A_\iota$ onto $B$, whose GNS representation is
faithful.
Let
$$ (A,\phi)=\freeprodi(A_\iota,\phi_\iota) $$
be the reduced amalgamated free product of C$^*$--algebras.
Then $A$ is exact.
\endproclaim
\demo{Proof}
$A$ is the inductive limit of C$^*$--algebras obtained by taking free products
of finite subfamilies of $\bigl((A_{\iota},\phi_\iota)\bigr)_{\iota\in I}$;
to see this, one can prove directly that the free product of a subfamily
embeds in the free product of the larger family by examining the Hilbert
C$^*$--modules on which these C$^*$--algebras act;
this also follows from the more general result
in~\cite{\BlanchardDykemaZZEmb}.
Since exactness is preserved under taking inductive limits, we may without loss
of generality assume that $I$ is finite.
Similarly, using Lemma~\SepSubalg{} and~\cite{\BlanchardDykemaZZEmb} as
described in Remark~\SepEmb, we may without loss of generality assume that $B$
and every $A_\iota$ is separable.

We shall use the same subspaces and completely positive maps as in
the proof of Lemma~\cpuinv, and the same notation.
Recall also that $P_\iota\oup$ is the projection from $E_\iota$ onto
$E_\iota\oup$.

It may happen that, $A_\iota\cap\Keu(E_\iota)\ne\{0\}$;
for technical reasons we would like to avoid this situation.
Let
$$ (D,\psi)=\freeprodi(A_\iota\otimes C(\Tcirc),\phi_\iota\otimes\tau), $$
where $C(\Tcirc)$ the C$^*$--algebra of all continuous functions on the circle
and where $\tau$ is the state on $C(\Tcirc)$ given by integration with respect
to Haar measure.
It is fairly straightforward in this example (and similar ones involving tensor
products) to show directly that $A$ is isomorphic to a C$^*$--subalgebra of
$D$;
moreover, this follows from the more general result found
in~\cite{\BlanchardDykemaZZEmb}.
Hence in order to show that $A$ is exact it will suffice to show that $D$ is
exact.
But each C$^*$--algebra $A_\iota\otimes C(\Tcirc)$ is exact and acting on
$L^2(A_\iota\otimes C(\Tcirc),\phi_\iota,\otimes\tau)$ contains no nonzero
compact operators, by Lemma~\NoCompacts.
Therefore, we may without loss of generality assume that
$A_\iota\cap\Keu(E_\iota)=\{0\}$ for all $\iota\in I$.

For $\iota\in I$ let $\Ah_\iota$ be the C$^*$--subalgebra of
$\Leu(E_\iota\oup)$ generated by
$\{\Po_\iota a\Po_\iota\mid a\in A_\iota\}\cup\Keu(E_\iota\oup)$.
We will now show that there is an exact sequence
$$ 0\to\Keu(E_\iota\oup)\to \Ah_\iota\to A_\iota\to 0 \tag{\KBAses} $$
with a completely positive unital splitting.
Let
$$ q_\iota\oup:\Leu(E_\iota\oup)\to\Leu(E_\iota\oup)/\Keu(E_\iota\oup) $$
be the quotient map.
Since $\Keu(E_\iota\oup)\subseteq \Ah_\iota$, we have the short exact sequence
$$ 0\to\Keu(E_\iota\oup)\to \Ah_\iota\overset{q_\iota\oup}\to\to
 q_\iota\oup(\Ah_\iota)\to0; $$
let us show that $q_\iota\oup(\Ah_\iota)\cong A_\iota$.
Clearly $q_\iota\oup(\Ah_\iota)$ is generated by
$\{q_\iota\oup(\Po_\iota a\Po_\iota)\mid a\in A_\iota\}$.
Since $\Po_\iota=1-\theta_{\xi_\iota,\xi_\iota}$,
$$ (\Po_\iota a_1\Po_\iota)(\Po_\iota a_2\Po_\iota)
-(\Po_\iota a_1a_2\Po_\iota) $$ 
is compact for all $a_1,a_2\in A_\iota$.
Therefore, the map
$$ A_\iota\ni a\mapsto q_\iota\oup(\Po_\iota a\Po_\iota) \tag{\Calkiniso} $$
is a $*$--homomorphism onto $q_\iota\oup(\Ah_\iota)$.
But if $a\in A_\iota$ and $\Po_\iota a\Po_\iota$ is compact then
$a\in\Keu(E_\iota)$ and hence $a=0$ by assumption.
Therefore the map in~(\Calkiniso) is an isomorphism from $A_\iota$ to
$q_\iota\oup(\Ah_\iota)$.
Now it is clear that the map
$$ A_\iota\ni a\mapsto\Po_\iota a\Po_\iota\in \Ah_\iota $$
is a completely positive unital splitting of the exact sequence~(\KBAses).

Fix $k\in\Nats$ and let us investigate $\Phi_k(A)$.
For $p\in\{1,2,\ldots,k\}$ and $\iota\in I$ consider the Hilbert $B$--modules
$$ \align
L_{p,\iota}&=\cases
\eta B&\text{if }p=k \\ \vspace{3ex}
\eta B\oplus\dsize\bigoplus\Sb
\ell\in\{1,2\ldots,k-p\} \\
\iota_1,\ldots,\iota_\ell\in I \\
\iota_1\ne\iota_2,\ldots,\iota_{\ell-1}\ne\iota_\ell \\
\iota_\ell\ne\iota \endSb
E_{\iota_1}\oup\otdts BE_{\iota_\ell}\oup &\text{if }p<k
\endcases \\ \vspace{2ex}
R_{p,\iota}&=\cases
\eta B&\text{if }p=1 \\ \vspace{3ex}
\dsize\bigoplus\Sb
\iota_1,\ldots,\iota_{p-1}\in I \\
\iota_1\ne\iota_2,\ldots,\iota_{p-2}\ne\iota_{p-1} \\
\iota_1\ne\iota \endSb
E_{\iota_1}\oup\otdts BE_{\iota_{p-1}}\oup &\text{if }p>1;
\endcases \endalign $$
As usual, $\eta B$ refers to a copy of the Hilbert $B$--module $B$ where the
identity element of $B$ is named $\eta$.
Let $V_{p,\iota}\in\Leu(L_{p,\iota}\otimes_BE_\iota\oup\otimes_BR_{p,\iota},E)$
be the isometry defined by erasing parenthesis and absorbing all occurrences of
$\eta$;
namely, given $e\in E_\iota\oup$, $\zeta_j\in E_{\iota_j}\oup$ and
$\zeta_j'\in E_{\iota_j'}\oup$ for
$\iota_1,\ldots,\iota_\ell,\iota_1',\ldots,\iota_{p-1}'\in I$ with
$\iota_j\ne\iota_{j+1}$, $\iota_j'\ne\iota_{j+1}'$, and using $\otimesdd$ for
the tensor product symbols in
$L_{p,\iota}\otimes_BE_\iota\oup\otimes_BR_{p,\iota}$, we let
$$ \alignat2
V_{p,\iota}:\,&\eta\otimesdd e\otimesdd\eta\mapsto e&\qquad&(\text{if }p=1) \\
&\eta\otimesdd e\otimesdd(\zeta_1'\otdt\zeta_{p-1}')
 \mapsto e\otimes\zeta_1'\otdt\zeta_{p-1}'&&(\text{if }p>1) \\
&(\zeta_1\otdt\zeta_\ell)\otimesdd e\otimesdd\eta
 \mapsto\zeta_1\otdt\zeta_\ell\otimes e&&(\text{if }p=1) \\
&(\zeta_1\otdt\zeta_\ell)\otimesdd e\otimesdd(\zeta_1'\otdt\zeta_{p-1}')
 \mapsto\zeta_1\otdt\zeta_\ell\otimes e
 \otimes\zeta_1'\otdt\zeta_{p-1}'&&(\text{if }p>1).
\endalignat $$
Note that for fixed $p$,
$\bigl(V_{p,\iota}(L_{p,\iota}\otimes_BE_\iota\oup\otimes_BR_{p,\iota})
V_{p,\iota}^*\bigr)_{\iota\in I}$ is a family of mutually orthogonal
complemented subspaces of $\Eto k$.

Let $a=a_1a_2\cdots a_q\in A$, where $a_j\in\Ao_{\iota_j}$,
$\iota_1\ne\iota_2,\ldots,\iota_{q-1}\ne\iota_q$.
Let $\zeta_j\in E_{\iota'_j}\oup$ ($1\le j\le n$) with
$\iota'_1\ne\iota'_2,\ldots,\iota'_{n-1}\ne\iota'_n$.
If $\iota_q\ne\iota'_1$ then
$$ a(\zeta_1\otdt\zeta_n)
=\widehat{a_1}\otdt\widehat{a_q}\otimes\zeta_1\otdt\zeta_n. $$
If $\iota_q=\iota'_1$ and $\iota_{q-1}\ne\iota'_2$ then
$$ \align
a(\zeta_1\otdt\zeta_n)
=\,&\widehat{a_1}\otdt\widehat{a_{q-1}}\otimes
 \langle\widehat{a_q^*},\zeta_1\rangle\zeta_2\otimes\zeta_3\otdt\zeta_n \\
&+\widehat{a_1}\otdt\widehat{a_{q-1}}\otimes
 (P_{\iota_q}\oup a_q\zeta_1)\otimes\zeta_2\otdt\zeta_n.
\endalign $$
Continuing in this way, and noting that, for example,
$$ \bigl\langle\widehat{a_{q-1}^*},
\langle\widehat{a_q^*},\zeta_1\rangle\zeta_2\bigr\rangle
=\langle\widehat{a_q^*}\otimes\widehat{a_{q-1}^*},
\zeta_1\otimes\zeta_2\rangle, $$
we see that if
$\iota_1=\iota'_1,\iota_{q-1}=\iota'_2,\ldots,\iota_{q-\ell+1}=\iota'_{\ell}$ but
$\iota_{q-\ell}\ne\iota'_{\ell+1}$, some $\ell<\min(q,n)$, then
$$ \align
a(\zeta_1\otdt\zeta_n)
=\,&\widehat{a_1}\otdt\widehat{a_{q-\ell}}\otimes
 \langle\widehat{a_q^*}\otdt\widehat{a_{q-\ell+1}^*},
 \zeta_1\otdt\zeta_\ell\rangle\zeta_{\ell+1}
 \otimes\zeta_{\ell+2}\otdt\zeta_n \\ \vspace{2ex}
&+\sum_{j=2}^\ell
 \topaligned \biggl(\widehat{a_1}&\otdt\widehat{a_{q-j}}\otimes \\
  &\otimes(P_{\iota_{q-j+1}}\oup a_{q-j+1}
   \langle\widehat{a_q^*}\otdt\widehat{a_{q-j+2}^*},
   \zeta_1\otdt\zeta_{j-1}\rangle\zeta_j)\otimes \\
  &\otimes\zeta_{j+1}\otdt\zeta_n\biggr) \endaligned \\ \vspace{2ex}
&+\widehat{a_1}\otdt\widehat{a_{q-1}}\otimes
 (P_{\iota_q}\oup a_q\zeta_1)\otimes\zeta_2\otdt\zeta_n.
\endalign $$
Making use of this and other related formulas and employing the convention
that a sum $\sum_{j=n}^mx_j$ is zero whenever $m<n$, we find that for
$a=a_1a_2\cdots a_q$ as above,
$$ \align
\Phi_k(a)&=\sum_{j=\max(0,q-k)}^{\min(q,k)}
 \theta_{\widehat{a_1}\otdt\widehat{a_j},\widehat{a_q^*}\otdt\widehat{a_{j+1}^*}}
 \tag{\PhikaTop} \\
 \vspace{2ex} 
&+\sum_{p=1}^{k-q}\topaligned
 &\left(\sum_{\iota\in I\backslash\{\iota_1\}}
  V_{p,\iota}\Bigl(\bigl(\theta_\eta1_{E_\iota\oup}
  \theta_{\widehat{a_q^*}\otdt\widehat{a_1^*}}^*\bigr)\otimes1_{R_{p,\iota}}
  \Bigr)V_{p,\iota}^*\right. \\
 &\quad\left.+\sum_{\iota\in I\backslash\{\iota_q\}}
  V_{p,\iota}\Bigl(\bigl(\theta_{\widehat{a_1}\otdt\widehat{a_q}}1_{E_\iota\oup}
  \theta_\eta^*\bigr)\otimes1_{R_{p,\iota}}
  \Bigr)V_{p,\iota}^*\right) \endtopaligned \\ \vspace{2ex}
&\!\!\!\!\!\!\!\!\!\!\!\!\!\!\!
 +\sum_{p=1}^{\min(k-1,k-\left[\frac{q+1}2\right])}
 \sum_{j=\max(1,p+q-k)}^{\min(q-1,k-p)}\,
 \sum_{\iota\in I\backslash\{\iota_j,\iota_{j+1}\}}
  V_{p,\iota}\Bigl(\bigl(\theta_{\widehat{a_1}\otdt\widehat{a_j}}1_{E_\iota\oup}
  \theta_{\widehat{a_q^*}\otdt\widehat{a_{j+1}^*}}^*\bigr)
  \otimes1_{R_{p,\iota}}\Bigr)V_{p,\iota}^* \\ \vspace{2ex}
&+\sum_{p=1}^{k-\left[\frac q2\right]}
 \sum_{j=\max(1,p+q-k)}^{\min(q,k-p+1)}
 V_{p,\iota_j}\Bigl(\bigl(\theta_{\widehat{a_1}\otdt\widehat{a_{j-1}}}
 P_{\iota_j}\oup a_jP_{\iota_j}\oup
 \theta_{\widehat{a_q^*}\otdt\widehat{a_{j+1}^*}}^*\bigr)\otimes1_{R_{p,\iota}}
 \Bigr)V_{p,\iota_j}^*, \tag{\PhikaBot}
\endalign $$
\vskip2ex
where the symbols
$$ \alignat2
\widehat{a_1}\otdt\widehat{a_j}&\text{ in line (\PhikaTop) }
 &\text{should be interpreted to mean }\xi&\text{ if }j=0, \\
\widehat{a_q^*}\otdt\widehat{a_{j+1}^*}&\text{ in line (\PhikaTop) }
 &\text{should be interpreted to mean }\xi&\text{ if }j=q, \\
\widehat{a_1}\otdt\widehat{a_{j-1}}&\text{ in line (\PhikaBot) }
 &\text{should be interpreted to mean }\eta&\text{ if }j=1, \\
\widehat{a_q^*}\otdt\widehat{a_{j+1}^*}&\text{ in line (\PhikaBot) }
 &\text{should be interpreted to mean }\eta&\text{ if }j=q.
\endalignat $$
Observe also that if $b\in B$ then
$$ \Phi_k(b)=\theta_{\widehat b,\xi}+\sum_{p=1}^k\sum_{\iota\in I}
V_{p,\iota}\bigl((\theta_\eta b\theta_\eta^*)
\otimes1_{R_{p,\iota}}\bigr)V_{p,\iota}^*. \tag{\Phikb} $$
The importance of these expressions is that they give
$\Phi_k(a_1\cdots a_q)$ and $\Phi_k(b)$ as sums of finitely many elements of
$$ \Keu(\Eto k)\cup\bigcup\Sb p\in\{1,\ldots,k\}\\ \iota\in I \endSb
V_{p,\iota}\bigl((\Keu(L_{p,\iota})\bowtie\Ah_\iota)
\otimes1_{R_{p,\iota}}\bigr)V_{p,\iota}^*. \tag{\KbtA} $$
The C$^*$--algebra generated by the set in~(\KbtA) therefore contains
$\Phi_k(A)$.
We will find a sort of composition series decomposition of this C$^*$--algebra
and thereby show that it is exact.

Given $p\in\{1,\ldots,k\}$ let
$$ D_p=\clspan\left(
\bigcup_{\iota\in I}V_{p,\iota}\Bigl(
\bigl(\Keu(L_{p,\iota})\bowtie\Ah_\iota\bigr)
\otimes1_{R_{p,\iota}}\Bigr)V_{p,\iota}\right). $$
Then $D_p$ is a C$^*$--subalgebra of $\Leu(\Eto k)$ that is isomorphic to 
$$ \bigoplus_{\iota\in I}\Keu(L_{p,\iota})\bowtie\Ah_\iota. $$
Moreover, it is not difficult to see that
$$ D_{p_1}D_{p_2}\subseteq D_{\min(p_1,p_2)}. \tag{\DpDp} $$
Let $\Ic_0=\Keu(\Eto k)$ and for every $p\in\{1,\ldots,k\}$ let
$\Ic_p=\Ic_{p-1}+D_p$.
Clearly~(\DpDp) implies that $\Ic_{p-1}$ is an ideal of $\Ic_p$.
We will show by induction on $p\in\{0,1,\ldots,k\}$ that $\Ic_p$ is a
C$^*$--algebra and that for every $p\ge1$ the short exact sequence
$$ 0\to\Ic_{p-1}\to\Ic_p\to\Ic_p/\Ic_{p-1}\to0 \tag{\IIses} $$
has a completely positive contractive splitting and
$$ \Ic_p/\Ic_{p-1}\cong\cases
\text{a C$^*$--subalgebra of }\bigoplus_{\iota\in I}\Keu\otimes A_\iota
 &\text{if }p<k \\
\bigoplus_{\iota\in I}A_\iota&\text{if }p=k.
\endcases \tag{\IoIiso} $$
Elementary theory of rings implies that the quotient $*$--algebra
$\Ic_p/\Ic_{p-1}$ is isomorphic to $D_p/(D_p\cap\Ic_{p-1})$, and the
isomorphism is seen to preserve the quotient norms.
Since $D_p/(D_p\cap\Ic_{p-1})$ is a C$^*$--algebra, it follows that $\Ic_p$ is
a C$^*$--algebra.

We have the commuting diagram
$$ \matrix
0\to&\Ic_{p-1}&\to&\Ic_p&\to&\Ic_p/\Ic_{p-1}&\to0 \\ \vspace{1.5ex}
&\cup&&\cup&&\uparrow \\ \vspace{1.5ex}
0\to&D_p\cap\Ic_{p-1}&\to&D_p&\to&D_p/(D_p\cap\Ic_{p-1})&\to0.
\endmatrix $$
In order to find a completely positive contractive splitting for the sequence
in the top row, it will suffice to find one for the sequence in the bottom row.
Since the projection $V_{p,\iota}V_{p,\iota}^*$ commutes with $D_p$ for each
$\iota\in I$, it will suffice to find a completely positive contractive
splitting for the exact sequence
$$ 0\to V_{p,\iota}^*(D_p\cap\Ic_{p-1})V_{p,\iota}
\to V_{p,\iota}^*D_pV_{p,\iota}\to
(V_{p,\iota}^*D_pV_{p,\iota})/(V_{p,\iota}^*(D_p\cap\Ic_{p-1})V_{p,\iota})
\to 0. \tag{\VDpVseq} $$
But
$$ V_{p,\iota}^*D_pV_{p,\iota}=\bigl(\Keu(L_{p,\iota})\bowtie\Ah_\iota\bigr)
\otimes1_{R_{p,\iota}} $$
while
$$ \Keu(L_{p,\iota}\otimes_BE_\iota\oup)\otimes1_{R_{p,\iota}}
\subseteq V_{p,\iota}^*(D_p\cap\Ic_{p-1})V_{p,\iota}
\subseteq\Keu(L_{p,\iota}\otimes_BE_\iota\oup)\bowtie\Leu(R_{p,\iota}). $$
Clearly
$$ \Keu(L_{p,\iota}\otimes_BE_\iota\oup)\bowtie\Leu(R_{p,\iota})
=\Keu(L_{p,\iota})\bowtie
\bigl(\Keu(E_\iota\oup)\bowtie\Leu(R_{p,\iota})\bigr) $$
and from Lemma~\capCompacts{} we have
$$ \bigl(\Keu(E_\iota\oup)\bowtie\Leu(R_{p,\iota})\bigr)
\cap\bigl(\Ah_\iota\otimes1_{R_{p,\iota}}\bigr)
=\Keu(E_\iota\oup)\otimes1_{R_{p,\iota}}; $$
hence
$$ V_{p,\iota}^*(D_p\cap\Ic_{p-1})V_{p,\iota}
=\bigl(\Keu(L_{p,\iota})\bowtie\Keu(E_\iota\oup)\bigr)\otimes1_{R_{p,\iota}}. $$
Therefore, using Lemma~\FoldedIdeals{} we see that $D_p/(D_p\cap\Ic_{p-1})$ is
isomorphic to a C$^*$--subalgebra of $\Keu\otimes A_\iota$, and, in light of
the completely positive contractive splitting we found for the short exact
sequence~(\KBAses), the short exact sequence~(\VDpVseq) has a completely
positive contractive splitting.
Putting these together, we obtain the isomorphisms~(\IoIiso) and a completely
positive contractive splitting for~(\IIses).

We can now finish the proof of the theorem.
By Lemma~\Kexact, the C$^*$--algebra $\Ic_0$ is exact.
The isomorphisms~(\IoIiso) imply that the quotients $\Ic_{p-1}/\Ic_p$ are exact
C$^*$--algebras.
Using the completely positively contractively split exact sequences~(\IIses)
and Lemma~\SplitExact, we can prove by induction that $\Ic_p$ is exact for all
$p\in\{0,1,\ldots,k\}$. 
But since $\Phi_k(A)$ is in the C$^*$--algebra generated by~(\KbtA), we see
that $\Phi_k(A)\subseteq\Ic_k$.
Therefore, we may use Lemma~\FactorHilmod{} together with the completely
positive unital maps $\Phi_k$ and $\Psi_k$ found in Lemma~\cpuinv{} to conclude
that $A$ is exact.
\QED

Recall that a discrete group is exact if and only if its reduced group
C$^*$--algebra is exact.
The following corollary is a result of Theorem~\exactfp{} and the fact,
expressed in equation~(\GpAmalgFP) in the introduction, that certain reduced
amalgamated free products of group C$^*$--algebras correspond to amalgamated
free products of groups.
\proclaim{Corollary \ExactGps}
Suppose $H$ is a group, $I$ is a set and for every $\iota\in I$, $G_\iota$ is
an exact group taken with the discrete topology and containing a copy of $H$ as
a subgroup.
Let $G=\fpiamalg HG_\iota$ be the free product of groups with amalgamation over
$H$.
Then $G$ is an exact group.
\endproclaim

\enddocument